\documentclass[11pt]{amsart}
\usepackage{amssymb, latexsym}
\theoremstyle{plain}
\newtheorem{theorem}{Theorem}

\newtheorem {lemma}{Lemma}

\theoremstyle{remark}

\newtheorem*{Remark 1}{Remark 1}
\newtheorem*{Remark 2}{Remark 2}
\newtheorem*{Remark 3}{Remark 3}
\newtheorem*{Remark 4}{Remark 4}

\numberwithin{equation}{section}

\begin{document}

\title[ Universal bound for  symmetric diffusion   with boundary flux]
 {Universal bound independent of geometry for solution to symmetric diffusion equation in exterior domain with boundary flux}

\author{Ross G. Pinsky}
\address{Department of Mathematics\\
Technion---Israel Institute of Technology\\
Haifa, 32000\\ Israel}
\email{ pinsky@math.technion.ac.il}
\urladdr{http://www.math.technion.ac.il/~pinsky/}

\date{}

\begin{abstract}
Fix $R>0$ and let
$B_R$ denote the ball of radius $R$ centered at the origin in $R^d$, $d\ge2$. Let $D\subset B_R$  be an open set with smooth boundary
and such that $R^d-\bar D$ is connected,
 and let 
$$
L=\sum_{i,j=1}^da_{i,j}\frac{\partial^2}{\partial x_i\partial x_j}-\sum_{i=1}^db_i\frac{\partial}{\partial x_i}
$$
be a second order elliptic operator. Consider 
the following linear heat equation in the exterior domain $R^d-\bar D$ with boundary flux:
\begin{equation*}
\begin{aligned}
&L u=0 \ \text{in}\ R^d-\bar D;\\
&a\nabla u\cdot \bar n=-h\ \text{on}\ \partial D;\\
&u>0 \ \text{is minimal},
\end{aligned}
\end{equation*}
where $h\gneq0$ is continuous, and where $\bar n$ is the unit inward normal to the domain $R^d-\bar D$.
The operator $L$ must possess a Green's function  in order that a solution  $u$ exist.
An important feature of the equation is that there is no a priori bound on the supremum $\sup_{x\in R^d-\bar D}u(x)$ of the solution
 exclusively in terms
of the boundary flux  $h$, the hyper-surface measure of $\partial D$ and the coefficients of $L$; rather the geometry
of $D\subset B_R$ plays an essential role.
However, we  prove that in the case that $L$ is a \it symmetric operator\rm\
with respect to some reference measure, then \it outside of\rm\ $\bar B_R$, the solution to \eqref{LHE} is uniformly bounded,
 independent of the particular choice of $D\subset B_R$. The proof uses a combination of analytic and probabilistic techniques.
\end{abstract}

\maketitle

\section{Introduction and Statement of  Results}
Fix $R>0$ once and for all and let
$B_R$ denote the ball of radius $R$ centered at the origin in $R^d$, $d\ge2$. Let $D\subset B_R$  be an open set with smooth boundary
and such that $R^d-\bar D$ is connected.
Let
\begin{equation}\label{nondivL}
L=\sum_{i,j=1}^da_{i,j}\frac{\partial^2}{\partial x_i\partial x_j}-\sum_{i=1}^db_i\frac{\partial}{\partial x_i}
\end{equation}
be a  strictly elliptic operator in $R^d-D$ with smooth coefficients $a=\{a_{i,j}\}_{i,j=1}^d$ and $b=\{b_i\}_{i=1}^d$.
Consider   the following linear heat equation in the exterior domain $R^d-\bar D$ with boundary flux:
\begin{equation}\label{LHE}
\begin{aligned}
&L u=0 \ \text{in}\ R^d-\bar D;\\
&a\nabla u\cdot \bar n=-h\ \text{on}\ \partial D;\\
&u>0 \ \text{is minimal},
\end{aligned}
\end{equation}
where $h\gneq0$ is continuous, and where $\bar n$ is the unit inward normal to the domain $R^d-\bar D$.
By minimal, we mean that  the solution  $u$ satisfies $u=\lim_{n\to\infty} u_n$, where
for $n>R$, $u_n$ solves
\begin{equation}\label{LHE-n}
\begin{aligned}
&L u=0 \ \text{in}\ B_n-\bar D;\\
&a\nabla u\cdot\bar n=-h\ \text{on}\ \partial D;\\
&u=0\ \text{on}\ \partial B_n.
\end{aligned}
\end{equation}
\medskip

We note   that a certain restriction must be placed on the operator $L$ in order that a solution  $u$ exist to \eqref{LHE}.
This will be discussed below.
An important feature of \eqref{LHE} is that there is no a priori bound on the supremum $\sup_{x\in R^d-\bar D}u(x)$ of the solution
 exclusively in terms
of the boundary flux  $h$, the hyper-surface measure of $\partial D$ and the coefficients of $L$; rather the geometry
of $D\subset B_R$ plays an essential role.
Here is a simple example.
\medskip

\bf \noindent Spherical Shell Example.\rm\
Consider a spherical shell centered at the origin in $R^3$ with the radii of the inner and outer boundary spheres given by $R_1$  and
$R_2$ respectively, where $R_1<R_2<R$. Now fix some direction $\theta_0\in S^2$ and puncture the shell with a spherical bullet of
radius $\frac1n$ centered in the $\theta_0$ direction. Thus, a cylindrical-like region of radius $\frac1n$ and length $R_2-R_1$ has been removed from the spherical shell.
Denote by $D_n$  the open set obtained by taking  this punctured spherical shell and deleting its boundary.
(If one insists that $D_n$ have a smooth boundary, one can smooth out the edges where the bullet enters and exits.)
Let $\Gamma_n= \partial B_{R_1}\cap (\bar D_n)^c$ denote the punctured part
of the inner sphere.  By the maximum principle, the solution $u_n$ to \eqref{LHE} with $L=\Delta$ and $h=1$ satisfies
$u_n(x)\ge v_n(x)$, for $|x|<R_1$, where $v_n$ is the solution of
\begin{equation*}
\begin{aligned}
&\Delta v_n=0 \ \text{in}\ B_{R_1};\\
&\nabla v_n\cdot \bar n=-1\ \text{on}\ \partial B_{R_1}-\bar\Gamma_n; \\
&v_n=0\ \text{on}\ \Gamma_n,
\end{aligned}
\end{equation*}
with $\bar n$ being   the inward unit normal to $B_{R_1}$ at $\partial B_{R_1}$.
By the maximum principle, $v_n(x)$ is increasing in $n$; let $v(x)\equiv\lim_{n\to\infty}v_n(x)$.
By Harnack's inequality, $v$ is either  finite everywhere or infinite everywhere in $B_{R_1}$.
If $v$ is finite, then $v$ must satisfy $\Delta v=0$  in $B_{R_1}$ and satisfy the Neumann boundary condition
$\nabla v\cdot \bar n=-1$ on all of $\partial B_{R_1}$. But this is impossible because solvability of the above equation
with  Neumann boundary data
requires the compatibility condition $\int_{\partial B_{R_1}}(\nabla v\cdot \bar n)d\sigma(x)=
-\int_{B_{R_1}}\Delta v\thinspace dx=0$, where $d\sigma$ denotes Lebesgue hyper-surface measure on $\partial B_{R_1}$. Thus, $v\equiv\infty$ and consequently,
$\lim_{n\to\infty}u_n(x)=\infty$, for $x\in B_{R_1}$.
Yet, trivially, the sequence  $\{|\partial D_n|\}_{n=1}^\infty$ is  bounded in $n$.

\medskip

However, we will prove that in the case that $L$ is a  symmetric operator
with respect to some reference measure, then for any $\delta>0$,  \it outside of\rm\ $B_{R+\delta}$, the solution to \eqref{LHE} is uniformly bounded,
 independent of the particular choice of $D\subset B_R$. More precisely, we will show that for any $x\in R^d-\bar B_R$, the solution $u(x)$ is \it bounded uniformly
over all  $D\subset B_R$\rm, in terms of

\noindent i. the boundary flux  $h$;

\noindent ii. the  hyper-surface measure of  $\partial D$ with respect
to the measure whose density with respect to Lebesgue measure is 
the trace of the above mentioned reference measure;

\noindent iii.   the behavior of the coefficients of $L$ outside of $\bar B_R$.

In particular, 
letting $\text{dist}_{\text{geod};D_n}(x_n,\bar B_{R+\delta})$ denote the length of a shortest path in $R^d-D_n$ from
$x_n$ to $\bar B_{R+\delta}$,
we will give an example of  a sequence of open sets $\{D_n\}_{n=1}^\infty\subset B_R$
and points $x_n\in B_R-D_n$ with\newline $\lim_{\delta\to0}\lim_{n\to\infty}\text{dist}_{\text{geod};D_n}(x_n,\bar B_{R+\delta})$
equal to  an arbitrarily small positive number and such that the solution $u_{D_n}$ to \eqref{LHE} with $D=D_n$ satisfies $\lim_{n\to\infty}u_{D_n}(x_n)=\infty$, yet
for all $x\in R^d-\bar B_R$, one has $\sup_nu_{D_n}(x)<\infty$.

The proof of our result will use a combination of analysis and probability. The analysis will consist
of representing the solution in terms of an appropriate Green's function and using symmetry, while the probability will consist of representing
this Green's function stochastically in terms of the occupation time of a diffusion process.

Before stating our result, we discuss the existence and uniqueness of a solution to \eqref{LHE}.
It is standard that the solution  $u_n$ to the linear equation \eqref{LHE-n} with co-normal boundary data at $\partial D$ and homogeneous Dirichlet data at $\partial B$
exists and is unique.
By the maximum principle, $u_n$ is positive off of $\partial B_n$ and attains its maximum on $\partial D$. The maximum principle also shows that
$u_n$ is nondecreasing in $n$. Thus, if $\lim_{n\to\infty}\sup_{x\in \partial D}u_n(x)<\infty$, then we obtain existence and uniqueness for the solution $u$ to \eqref{LHE}.
The maximum principle shows that $u$ attains its maximum
on the boundary $\partial D$.

With regard to existence, we begin with a physical description.
When $d=3$, $u$ can be thought of as the equilibrium quantity of a reactant  after having  undergone a long period of $L$-diffusion
and convection in an exterior domain which is being supplied with the reactant via a boundary  flux $h$,  and where  complete and instantaneous   absorption
occurs far away. Note that in \eqref{nondivL}, the drift term $b$ has been written with a minus sign;
with this convention, the reactant is being convected in the direction $b$.
If this convection vector field is pointing away from $D$ with sufficient strength, then
 the convection of the reactant from the boundary $\partial D$ into the domain $R^3-\bar D$ will
not allow for an equilibrium state.

In order to obtain existence,  in fact it  is necessary and sufficient to assume that the operator $L$  possesses an appropriate
Green's function.

\noindent \bf Assumption G.\rm\ The operator $L$ with the co-normal  boundary condition at $\partial D$ possesses a
Green's function.
\medskip

Recall that the Green's function $G^D_{\text{Neu}}(x,y)$ for an operator  $L$ with the co-normal boundary condition at $\partial D$
is the minimal positive function $g(x,y)$ satisfying the following conditions: for each $y\in R^d-\bar D$ the function $g(\cdot, y)$ satisfies
  $Lg(\cdot, y)=-\delta_y$ in $R^d-\bar D$ and satisfies the homogeneous co-normal boundary condition at $\partial D$:
  $a\nabla g(\cdot,y)\cdot\bar n=0$ on $\partial D$.
The Green's function for an operator $L$ on all of $R^d$ is the minimal positive function $g(x,y)$ satisfying
$Lg(\cdot, y)=-\delta_y$ in $R^d$, for each $y\in R^d$.
We make several remarks concerning the Green's function.

\bf\noindent Remark 1.\rm\ Assumption G is equivalent to each of the following assumptions:

\noindent \it i.\rm\ If  $L$ is extended to be a smooth strictly elliptic operator on all of $R^d$, then this extension possesses a Green's function;

\noindent \it ii.\rm\ For $n>R$, let $V_n$ denote the solution to the equation
\begin{equation}\label{V1}
\begin{aligned}
&LV_n=0\ \text{in}\ B_n-\bar B_R;\\
& V_n=1\ \text{on}\ \partial B_R,\  V_n=0\ \text{on}\ \partial B_n.
\end{aligned}
\end{equation}
Then
\begin{equation}\label{V2}
V\equiv\lim_{n\to\infty}V_n
\end{equation}
is not the constant function 1.

\noindent \it iii.\rm\ The diffusion process $X(t)$ in $R^d-D$ corresponding to the operator
$L$ and with co-normal reflection at $\partial D$ is transient, that is,
$P_x(\lim_{t\to\infty}|X(t)|=\infty)=1$, or equivalently,
$P_x(\tau_R<\infty)<1$, for $x\in R^d-\bar B_R$, where
$P_x$ denotes probabilities for the diffusion starting from $x$
and $\tau_R=\inf\{t\ge0:X(t)\in\bar B_R\}$ is the first hitting time of the ball $\bar B_R$.
 Indeed, the function $V$ defined
in \eqref{V1}-\eqref{V2} satisfies 
\begin{equation}\label{Vprob}
V(x)=P_x(\tau_R<\infty).
\end{equation}
For details, see \cite{P}.

\medskip

\noindent \bf Remark 2.\rm\ If Assumption G is not in force, then the solution $u_n$ to  \eqref{LHE-n} will satisfy
$\lim_{n\to\infty}u_n=\infty$.

\noindent \bf Remark 3.\rm\  For generic choices of $b$ such that $L$  satisfies Assumption G,  the solution $u$ to
\eqref{LHE}  will satisfy  $\lim_{x\to\infty}u(x)=0$.
Necessarily, one has $\liminf_{x\to\infty}u(x)=0$.  One can
obtain $\limsup_{x\to\infty}u(x)>0$ by choosing the convection term $b$   pointing very strongly away from the origin when $x$ is in certain  sectors.
(The restriction to certain sectors is necessary for otherwise the Green's function would not exist.)
\medskip

Our result requires that $L$ be symmetric with respect to a weight $e^Q$; that is,
$L$ must be of the form
\begin{equation}\label{Lsym}
L=e^{-Q}\nabla\cdot e^{Q}a\nabla=\nabla\cdot a\nabla+a\nabla Q\nabla.
\end{equation}
(Note that such an $L$ is of the form \eqref{nondivL} with $b_i=-\sum_{j=1}^d\frac{\partial a_{i,j}}{\partial x_j}-\sum_{j=1}^da_{i,j}\frac{\partial Q}{\partial x_j}$.)
Note that the boundary condition in \eqref{nondivL} involves the co-normal derivative;
as is well-known, the operator  $L$ in \eqref{Lsym} with
the homogeneous co-normal boundary condition is symmetric with resect to the weight $e^Q$:
 $\int_{R^d-D}gLfe^Qdx=\int_{R^d-D}fLge^Qdx$, for all smooth, compactly supported $f,g$ that satisfy
 $a\nabla f\cdot\bar n=a\nabla g\cdot \bar n=0$ on $\partial D$.

We first present a theorem for the case $L=\Delta$ with $d\ge3$.
(We must restrict to $d\ge3$ because there is no Green's function when $d=2$.)
We can
obtain tighter results in this case  than in the case of generic $L$. We note that the Green's
function for $\Delta$ in $R^d$ is $G(x,y)=\frac{|x-y|^{2-d}}{(d-2)\omega_d}$ where
$\omega_d$ denotes the Lebesgue hyper-surface measure of the unit sphere in $R^d$.
 Lebesgue hyper-surface measure on $\partial D$ is denoted
 by $d\sigma$.

\begin{theorem}\label{LHE-Lap}
Let $L=\Delta$ in $R^d$, $d\ge 3$. Let $R>0$ and assume that $D\subset B_R$. Let $\omega_d$ denote the surface
measure of the unit sphere in $R^d$.
Then for every $\gamma>1$, the solution $u$ to \eqref{LHE} satisfies
\begin{equation*}
\left(\int_{\partial D}h(z)d\sigma(z)\right)
 \frac{c^-_{\gamma,d}\thinspace|x|^{2-d}}{(d-2)\omega_d}\le u(x)\le
\left(\int_{\partial D}h(z)d\sigma(z)\right)
 \frac{c^+_{\gamma,d}\thinspace|x|^{2-d}}{(d-2)\omega_d},
\end{equation*}
for $|x|\ge \gamma R$,
where $c_{\gamma,d}^\pm$ are independent of $D$ and $R$ and satisfy  $\lim_{\gamma\to\infty}c^\pm_{\gamma,d}=1$.
\end{theorem}
For the case of a generic operator $L$, we need one  more definition before we can state the result.
Let $G^R_{\text{Dir}}$ denote the Green's function for $L$ in $R^d-\bar B_R$ with the Dirichlet boundary condition at $\partial B_R$.
That is, $G^R_{\text{Dir}}(x,y)$ is the minimal positive function $g(x,y)$ satisfying
  $Lg(\cdot, y)=-\delta_y$ in $R^d-\bar B_R$ and satisfying  the zero  Dirichlet boundary condition at $\partial B_R$,
for each $y\in R^d-\bar B_R$. As an aside, we note that the Green's function $G^R_{\text{Dir}}$ always exists, even if Assumption G is not
satisfied.
\begin{theorem}\label{LHE-gen}
Assume that the operator $L$ satisfies Assumption G and is symmetric with respect to the weight function $e^Q$ as in \eqref{Lsym}.
Let $R>0$ and assume that $D\subset B_R$. Let $R'>R$. Then the solution $u$ to \eqref{LHE} satisfies
\begin{equation}\label{twosided}
\begin{aligned}
&\left(\int_{\partial D}h(y)e^{Q(y)}d\sigma(y)\right)\frac{\min_{|z|=R'}G^R_{\text{Dir}}(z,x)}{1-\min_{|z|=R'}V(z)}e^{-Q(x)}
\le u(x)\le \\
&\left(\int_{\partial D}h(y)e^{Q(y)}d\sigma(y)\right)\frac{\max_{|z|=R'}G^R_{\text{Dir}}(z,x)}{1-\max_{|z|=R'}V(z)}e^{-Q(x)},\
\text{for}\ |x|>R',
\end{aligned}
\end{equation}
where $G^R_{\text{Dir}}$ is the Green's function for $L$ in $R^d-B_R$ with the  Dirichlet boundary condition at $\partial B_R$,
and $V$ is as in \eqref{V1}-\eqref{Vprob}.
\end{theorem}
\noindent\bf Remark 1.\rm\ Recall that $V\not\equiv1$ is equivalent to Assumption G.
By the maximum principal, $V<1$ in $R^d-\bar B_R$. In the generic
case one has $\lim_{x\to\infty}V(x)=0$. Necessarily one has $\liminf_{x\to\infty}V(x)=0$. One can obtain
$\limsup_{x\to\infty}V(x)>0$ by choosing $Q$ in the manner noted in Remark 3 following Assumption G.
If $\lim_{x\to\infty}V(x)=0$, then of course $\lim_{R'\to\infty}\frac{1-\min_{|z|=R'}V(z)}{1-\max_{|z|=R'}V(z)}=1$.
For certain classes of operators one has
 $\lim_{x\to\infty}\frac{\min_{|z|=R'}G^R_{\text{Dir}}(z,x)}{\max_{|z|=R'}G^R_{\text{Dir}}(z,x)}=1$.
If the above two limits hold, then by  choosing $R'$ large, the ratio of the left hand side
to the right hand side of \eqref{twosided} can be made arbitrarily close to 1 for large $|x|$.

\noindent \bf Remark 2.\rm\
One can of course choose a sequence $\{D_n\}_{n=1}^\infty$ of domains satisfying
$D_n\subset B_R$, for all $n$, and $\lim_{n\to\infty}|\partial D_n|=\infty$.
Letting $u_n$ denote the solution to \eqref{LHE} with, say, $h=1$ on $\partial D_n$,
it follows  that $\lim_{n\to\infty}u_n(x)=\infty$,
for  $|x|> R$.

Theorems \ref{LHE-Lap} and \ref{LHE-gen}
show that given an $M>0$,  then for any $\delta>0$, the solution to \eqref{LHE} is bounded in $R^d- B_{R+\delta}$ uniformly over all
$D\subset B_R$ satisfying $|\partial D|\le M$ and over all fluxes $h$ satisfying $|h|\le M$, even though
the solution can be arbitrarily large at points inside $B_R$ and very
close to $\partial B_R$.
For example,  consider the  case of the punctured spherical shell, defined above, except  let the inner and outer radii, $R_1$ and $R_2$, also
 depend on $n$, setting $R_1=R-\frac2n$ and
 $R_2=R-\frac1n$.
 Then for any $\delta>0$, the  solutions $\{u_{D_n}\}_{n=1}^\infty$ will be uniformly bounded over $R^d-B_{R+\delta}$,
 even though for  $x_n$  satisfying $|x_n|=1-\frac3n$ and $\text{arg}(\frac{x_n}{|x_n|})=\theta\neq\theta_0$, one has that  $\lim_{n\to\infty}u_{D_n}(x_n)=\infty$.
Recalling the definition of $\text{dist}_{\text{geod};D_n}(\cdot,\cdot)$, note
 that $\lim_{\delta\to0}\lim_{n\to\infty}\text{dist}_{\text{geod};D_n}(x_n,\bar B_{R+\delta})$ can be made as small as one likes
by choosing $\theta$ sufficiently close to $\theta_0$.

The proofs of Theorems \ref{LHE-Lap} and \ref{LHE-gen} depend in an essential way on the fact that the operator $L$
is symmetric.
Of course, the statement of the result also depends on symmetry, as the weight function $e^Q$ appears in the boundary integrals.

\noindent \bf Open Question.\rm\ In the case that the operator $L$ is not symmetric, is the size of the solution to \eqref{LHE} at  $x\in R^d-\bar B_R$ governed independently of the geometry of $D$? Or alternatively, can one, say,  give an example of a class
of domains $\{D_n\}_{n=1}^\infty \subset B_R$ with $\sup_n|\partial D_n|<\infty$ and such that the corresponding solutions
 $\{u_{D_n}\}_{n=1}^\infty$ to \eqref{LHE} with, say, $h\equiv1$ satisfy  $\sup_n\sup_{|x|>R+\delta}u_{D_n}(x)=\infty$, for some $\delta>0$?
\medskip

The proof of Theorem \ref{LHE-Lap} is given in section \ref{sec:LHE-Lap}. In section \ref{sec:gen} we sketch how to amend the
 proof of Theorem \ref{LHE-Lap} to obtain the proof of Theorem \ref{LHE-gen}.

\section{Proof of Theorem \ref{LHE-Lap}}\label{sec:LHE-Lap}
\begin{proof}
For $n>R$, let $u_n$ denote the solution to
\begin{equation}\label{PN}
\begin{aligned}
&\Delta u_n=0 \ \text{in}\ B_n-\bar D;\\
&\nabla u_n\cdot n=-h\ \text{on}\ \partial D;\\
&u_n=0\ \text{on}\ \partial B_n.
\end{aligned}
\end{equation}
Let $G^D_n(x,y)$ denote the Green's function  for  $\Delta$
in $B_n-\bar D$
with the  Dirichlet boundary condition
at $\partial B_n$ and the  Neumann boundary condition at $\partial D$.
One has
\begin{equation}
\begin{aligned}
&\Delta_x G^D_n(x,y)=-\delta_y\ \text{for}\ x,y\in B_n-\bar D;\\
&G^D_n(x,y)=0, \ \text{for}\ x\in \partial B_n, \ y\in B_n-\bar D;\\
&\frac{\partial G^D_n}{\partial \bar n_x}(x,y)=0,\ \text{for}\ x\in\partial D, \ y\in B_n-\bar D.
\end{aligned}
\end{equation}
 By the maximum principle, $G^D_n$ is increasing in $n$; let
 $G^D_{\text{Neu}}(x,y)=\lim_{n\to\infty}G^D_n(x,y)$.  This limiting function is the  Green's function
 for $\Delta$ in $R^d-\bar D$ with the Neumann boundary condition at $\partial D$.

Recall the fundamental property of the Green's function: if $v$ is a smooth function in $B_n-D$ and
continuous up to $\partial B_n$, then
$$
\begin{aligned}
&v(x)=-\int_{B_n-D}(\Delta v)(y)G_n^D(x,y)dy+\int_{\partial B_n}v(y)(\nabla_y G_n^D\cdot \bar n)(x,y)d\sigma(y)-\\
&\int_{\partial D}(\nabla v\cdot \bar n)(y)G_n^D(x,y)d\sigma(y),
\end{aligned}
$$
where $\bar n$ denotes the exterior unit normal vector to $D$ at $\partial D$ and the interior unit  normal
vector to $B_n$ at $\partial B_n$.
 Since $u_n$ is harmonic and since  $u_n$ vanishes on $\partial B_n$ and $\nabla u_n\cdot \bar n=-h$ on $\partial D$, one has the representation
$$
u_n(x)=\int_{\partial D}G^D_n(x,y)h(y)d\sigma(y), \ \text{for}\  x\in B_n-D.
$$
 Letting $n\to\infty$ gives
 \begin{equation}\label{formula}
 u(x)=\int_{\partial D}G^D_{\text{Neu}}(x,y)h(y)d\sigma(y), \ \text{for}\  x\in R^d-\bar D.
 \end{equation}
 To complete the proof of the theorem, we  appeal to spectral theory and then to the
   probabilistic representation
of  the Green's function $G^D_{\text{Neu}}(x,y)$.

Since the operator $\Delta$ in $B_n-\bar D$ with the Dirichlet
boundary condition at $\partial B_n$ and the Neumann boundary
condition
 at $\partial D$ is symmetric, it follows that the Green's function $G^D_n$, which is the integral kernel
 of the inverse operator, is symmetric; that is,
 $G^D_n(x,y)=G^D_n(y,x)$. Thus, also
 \begin{equation}\label{sym}
 G^D_{\text{Neu}}(x,y)=G^D_{\text{Neu}}(y,x).
\end{equation}

We now turn to the probabilistic representation of the Green's function.
 Let $B(t)$ be a  $d$-dimensional Brownian motion
 in $R^d-D$,  normally reflected at $\partial D$, and  corresponding to the operator $\Delta$.
(A standard Brownian motion $\beta(t)$ corresponds to the operator $\frac12\Delta$; our $B(t)$ can be obtained as
$\sqrt2\beta(t)$, or alternatively, as $\beta(2t)$.)
Let $P_x$ denote probabilities and let $E_x$ denote the corresponding  expectations for the Brownian motion starting from $x\in R^d-D$.
 For $x\in R^d-D$, define the expected occupation measure
 by $\mu_x^D(A)=E_x\int_0^\infty1_A(B(t))dt$, for Borel sets $A\subset R^d-D$.
 Since $d\ge3$, the Brownian motion is transient
 (that is, $P_x(\lim_{t\to\infty}|B(t)|=\infty)=1$) and from this one can show that   $\mu_x^D(A)<\infty$ for all bounded $A$.
 The measure $\mu_x^D(dy)$ possesses  a density
and the density is given by $G^D(x,y)$ \cite{P}.
From now on  we will write $G^D(x,A)\equiv \mu_x^D(A)$.
Using this probabilistic representation,  we will  show that for $\gamma>1$,
\begin{equation}\label{key}
\begin{aligned}
&G^D_{\text{Neu}}(x,y)\le  \frac1{(d-2)\omega_d}c^+_{\gamma,d}|y|^{2-d},
\ \text{for}\ x\in\partial D, |y|\ge \gamma R;\\
&G^D_{\text{Neu}}(x,y)\ge\frac1{(d-2)\omega_d}c^-_{\gamma,d}|y|^{2-d},
\ \text{for }\ x\in \partial D, |y|\ge \gamma R,
\end{aligned}
\end{equation}
where $\lim_{\gamma\to\infty} c^+_{\gamma,d}=\lim_{\gamma\to\infty}c^-_{\gamma,d}=1$,  $c^-_{\gamma,d}>0$
and $c^\pm_{\gamma,d}$ are independent of $D$ and $R$.
The theorem then follows from \eqref{formula}, \eqref{sym} and \eqref{key}.

To prove \eqref{key}, we define a sequence of hitting times for the Brownian motion.
Let $\gamma>\rho>1$.
Define  $\tau_1=\inf\{t\ge0:|B(t)|=\rho R\}$,  and then by induction define
$\tau_{2n}=\inf\{t>\tau_{2n-1}:|B(t)|=R\}$ and $\tau_{2n+1}=\inf\{t>\tau_{2n}: |B(t)|=\rho R\}$.
In words, $\tau_1$ is the first time the Brownian motion hits $\partial B_{\rho R}$, $\tau_2$ is the first time after
$\tau_1$ that the Brownian motion hits $\partial B_R$, $\tau_3$ is the first time after $\tau_2$ that the Brownian motion
hits $\partial B_{\rho R}$, etc. Since the Brownian motion is transient,
almost surely only a finite number of the $\tau_n$ will be finite.
For $x\in\partial D$ and  $A\subset R^d- B_{\gamma R}$, we have
\begin{equation}\label{circuits}
G^D_{\text{Neu}}(x,A)=E_x\int_0^\infty1_A(B(t))dt=\frac12\sum_{n=1}^\infty E_x\int_{\tau_{2n-1}}^{\tau_{2n}}1_A(B(t))dt.
\end{equation}
This equality holds because for times $s\in[\tau_{2n},\tau_{2n+1}]$ one has $|B(s)|\le \rho R$ and thus $B(s)\not\in A$.

Note that
  for  $z\in\partial B_{\rho R}$, one has $P_z(\tau_1=0)=1$ and thus under $P_z$ one has that $\tau_2$ is the first hitting time of $B_R$.
  Let $\phi(r)=(\frac Rr)^{d-2}$, $r=|x|$,  be the radially symmetric harmonic function in $R^d-\bar B_R$ which equals 1 on the boundary and decays to 0 at $\infty$.
It is well-known that starting from $z$ with $|z|>r$, the probability that the first hitting time of $\partial  B_R$ is finite
is $\phi(|z|)$; thus, $P_z(\tau_2<\infty)=\rho^{2-d}$, for $z\in\partial B_{\rho R}$.
By the strong Markov property, conditioned on $\tau_{2n-1}<\infty$, the probability that $\tau_{2n}<\infty$
is again $\rho^{2-d}$. Thus, one has
\begin{equation}\label{hitting}
P_x(\tau_{2n-1}<\infty)=\rho^{(2-d)(n-1)}, \ \text{for}\ x\in \partial D.
\end{equation}

Also from the strong Markov property, the conditional expectation\newline $E_x(\int_{\tau_{2n-1}}^{\tau_{2n}}1_A(B(t))dt|\tau_{2n-1}<\infty)$ is some averaging of the values
of \newline $E_z\int_0^{\tau_2}1_A(B(t))dt$, as $z$ varies over $\partial B_{\rho R}$. That is, there exists a probability measure
$\nu_{n,x}$ on $\partial B_{\rho R}$ such that
\begin{equation}\label{avg}
E_x(\int_{\tau_{2n-1}}^{\tau_{2n}}1_A(B(t))dt|\tau_{2n-1}<\infty)=\int_{\partial B_{\rho R}}\left(E_z\int_0^{\tau_2}1_A(B(t))dt\right)\nu_{n,x}(dz).
\end{equation}

Now let $W(t)$ be a   Brownian motion in all of $R^d$ corresponding to the operator $\Delta$ and let $\mathcal {E}_x$  denote the expectation for this
Brownian motion starting from $x$.
Since the  Brownian motion reflected at $\partial D$ and the  Brownian motion on all of $R^d$  behave the same  when they are in $R^d-B_R$, one has
\begin{equation}\label{2bm}
E_z\int_0^{\tau_2}1_A(B(t))dt=\mathcal{E}_z\int_0^{\tau_2}1_A(W(t))dt,\ \text{ for}\ z\in\partial B_{\rho R}.
\end{equation}

We now prove the upper bound  in  \eqref{key}.
From \eqref{2bm} we have
\begin{equation}\label{compare}
E_z\int_0^{\tau_2}1_A(B(t))dt\le \mathcal{E}_z\int_0^\infty1_A(W(t))dt, \ \text{ for}\  z\in\partial B_{\rho R}.
\end{equation}
Recall that the Green's function for $\Delta$ on all of $R^d$ is given by
$G(x,y)\equiv\frac{|x-y|^{2-d}}{(d-2)\omega_d}$, where $\omega_d$ denotes the surface measure of the unit sphere in $R^d$.
The probabilistic representation of the Green's function described above also holds for the  Brownian motion in all of $R^d$; that is,
$G(x,y)$ is the density of the measure $\mu_x(A)=\mathcal{E}_x\int_0^\infty1_A(W(t))dt$. Thus,
one has
\begin{equation}\label{free}
\mathcal{E}_z\int_0^\infty1_A(W(t))dt=\int_AG(z,y)dy=\int_A\frac{|z-y|^{2-d}}{(d-2)\omega_d}dy.
\end{equation}
Since $E_x(\int_{\tau_{2n-1}}^{\tau_{2n}}1_A(B(t))dt|\tau_{2n-1}=\infty)=0$, we have from
\eqref{circuits}-\eqref{free} that
\begin{equation}\label{final}
\begin{aligned}
&G^D_{\text{Neu}}(x,A)=E_x\int_0^\infty1_A(B(t))dt=\sum_{n=1}^\infty E_x\int_{\tau_{2n-1}}^{\tau_{2n}}1_A(B(t))dt=\\
&\sum_{n=1}^\infty E_x(\int_{\tau_{2n-1}}^{\tau_{2n}}1_A(B(t))dt|\tau_{2n-1}<\infty)P_x(\tau_{2n-1}<\infty)\le\\
&\left(\sup_{z\in \partial B_{\rho R}} \mathcal{E}_z\int_0^\infty1_A(W(t))dt\right)\sum_{n=1}^\infty P_x(\tau_{2n-1}<\infty)=\\
&\left(\sup_{z\in \partial B_{\rho R}}\int_A\frac{|z-y|^{2-d}}{(d-2)\omega_d}dy\right)\sum_{n=1}^\infty \rho^{(2-d)(n-1)}\le\\
&\int_A\frac1{(d-2)\omega_d}\frac1{1-\rho^{2-d}}\left(\frac\gamma{\gamma-\rho}    \right)^{d-2}|y|^{2-d}dy, \ \text{for}\ x\in\partial D,\ A\subset R^d-B_{\gamma R},
\end{aligned}
\end{equation}
where in the last inequality we have used the fact that
  $|y|\le \frac{\gamma}{\gamma-\rho}|y-z|$, for $|z|=\rho R$ and $|y|\ge\gamma R$.
Now  $\rho\in(1,\gamma)$ is a free parameter.
One can check that
 $$
 \sup_{\rho\in(1,\gamma)}(1-\rho^{2-d})(\gamma-\rho)^{d-2}=
 \frac{(\gamma^{\frac{d-2}{d-1}}-1)(\gamma-\gamma^{\frac1{d-1}})^{d-2}}{\gamma^{\frac{d-2}{d-1}}}.
$$
The maximum is attained at $\rho=\gamma^\frac1{d-1}$.
Using this value of $\rho$ in the right hand side of \eqref{final} gives
\begin{equation}\label{finalagain}
G^D_{\text{Neu}}(x,A)\le\int_A\frac{c^+_{\gamma,d}}{(d-2)\omega_d}|y|^{2-d}dy,
\end{equation}
for  $x\in\partial D$ and   $A\subset R^d-B_{\gamma R}$, where
$$
c^+_{\gamma,d}=\frac{\gamma^{\frac{d-2}{d-1}}\gamma^{d-2}}{(\gamma^{\frac{d-2}{d-1}}-1)(\gamma-\gamma^{\frac1{d-1}})^{d-2}}.
$$
Note that $\lim_{\gamma\to\infty}c^+_{\gamma,d}=1$.
Now the upper bound in \eqref{key}  follows from \eqref{finalagain}.

We now prove the lower bound  in \eqref{key}.
Let $G^R_{\text{Dir}}(x,y)$
 denote the Green's function for $\Delta$ in $R^d-\bar B_R$ with the Dirichlet
boundary condition at $\partial B_R$.
The probabilistic representation of the Green's function gives
\begin{equation}\label{Greenconn}
G^R_{\text{Dir}}(z,A)=\mathcal{E}_z\int_0^{\tau_2}1_A(W(t))dt, \ \text{for}\ |z|>R.
\end{equation}
Using reflection with respect to $\partial B_R$ allows one to calculate  the Green's function explicitly \cite{J}:
\begin{equation}\label{GreenR}
\begin{aligned}
&G^R_{\text{Dir}}(x,y)=\frac1{(d-2)\omega_d}|x-y|^{2-d}-\frac1{(d-2)\omega_d}(\frac{|y|}R)^{2-d}|x-\frac{R^2}{|y|^2}y|^{2-d}=\\
&\frac1{(d-2)\omega_d}|x-y|^{2-d}-\frac1{(d-2)\omega_d}\left|\frac{|y|}Rx-\frac R{|y|}y\right|^{2-d}.
\end{aligned}
\end{equation}
(In the proof of the upper bound, we could have used this Green's function and \eqref{2bm} instead of the Green's function
$G(x,y)$ for all of $R^d$ and \eqref{compare}, and this would have yielded a slightly smaller value of $c^+_{\gamma,d}$. However,
it was simpler to work with $G(x,y)$. For the lower bound we have no choice but to work with $G^R_{\text{Dir}}(x,y)$.)

For $|x|=\rho R$ and $|y|=\gamma' R$, with $\gamma'\ge \gamma$,
one has
$$
\left|\frac{|y|}Rx-\frac R{|y|}y\right|\ge \frac{|y||x|}R-R=R(\gamma'\rho-1)
$$
and
$$
|x-y|\le |x|+|y|=R(\gamma'+\rho).
$$
Using this with \eqref{GreenR}, one has
\begin{equation}\label{lowerbound}
\frac{(d-2)\omega_dG^R_{\text{Dir}}(x,y)}{|y|^{2-d}}\ge (\frac{\gamma'}{\gamma'+\rho})^{d-2}-(\frac{\gamma'}{\gamma'\rho-1})^{d-2},\text{for}\ |x|=\rho R,\  |y|=\gamma' R.
\end{equation}
One can check that the right hand side of \eqref{lowerbound} is increasing in $\gamma'$; thus
\begin{equation}\label{lowerboundagain}
\frac{(d-2)\omega_dG^R_{\text{Dir}}(x,y)}{|y|^{2-d}}\ge (\frac{\gamma}{\gamma+\rho})^{d-2}-(\frac{\gamma}{\gamma\rho-1})^{d-2},\ \text{for}\ |x|=\rho R,\  |y|\ge\gamma R.
\end{equation}

Then similar to \eqref{final}, but using \eqref{2bm} and \eqref{Greenconn} instead of \eqref{compare} and \eqref{free},  we have
\begin{equation}\label{finallower}
\begin{aligned}
&G^D_{\text{Neu}}(x,A)=E_x\int_0^\infty1_A(B(t))dt=\sum_{n=1}^\infty E_x\int_{\tau_{2n-1}}^{\tau_{2n}}1_A(B(t))dt=\\
&\sum_{n=1}^\infty E_x(\int_{\tau_{2n-1}}^{\tau_{2n}}1_A(B(t))dt|\tau_{2n-1}<\infty)P_x(\tau_{2n-1}<\infty)\ge\\
&\left(\inf_{z\in \partial B_{\rho R}} \mathcal{E}_z\int_0^{\tau_2}1_A(W(t))dt\right)\sum_{n=1}^\infty P_x(\tau_{2n-1}<\infty)=\\
&\left(\inf_{z\in \partial B_{\rho R}}\int_AG^R_{\text{Dir}}(z,y)dy\right)\sum_{n=1}^\infty \rho^{(2-d)(n-1)}=
\left(\inf_{z\in \partial B_{\rho R}}\int_AG^R_{\text{Dir}}(z,y)dy\right)\frac1{1-\rho^{2-d}}
\ge\\
&\int_A\frac1{(d-2)\omega_d}\frac1{1-\rho^{2-d}}
\left((\frac{\gamma}{\gamma+\rho})^{d-2}-(\frac{\gamma}{\gamma\rho-1})^{d-2}\right)|y|^{2-d}dy, \\
& \text{for}\ x\in\partial D,\ A\subset R^d-B_{\gamma R},
\end{aligned}
\end{equation}
where the inequality follows from  \eqref{lowerboundagain}.

As a function of $\rho\in(1,\gamma)$, the expression $\frac1{1-\rho^{2-d}}
\left((\frac{\gamma}{\gamma+\rho})^{d-2}-(\frac{\gamma}{\gamma\rho-1})^{d-2}\right)$ is non-positive if $\gamma\le1+\sqrt2$.
For $\gamma>1+\sqrt2$, this expression has a positive maximum. The formula for the maximum is complicated
and doesn't add much so we will just define $c_{\gamma,d}^-$ as follows.
For $\gamma$ sufficiently large so that the above expression is positive
for $\rho=\gamma^{\frac1 {d-1}}$, let this value be $c_{\gamma,d}^-$.
 Then, in particular, for some $\gamma_0>1$,  one has
 $$
c_{\gamma,d}^-=\frac1{1-\gamma^{\frac{2-d}{d-1}}}
\left(\left(\frac\gamma{\gamma+\gamma^{\frac1{d-1}}}\right)^{d-2}-\left(\frac\gamma{\gamma^{\frac d{d-1}}-1}\right)^{d-2}\right),
\  \text{for}\  \gamma\ge \gamma_0.
 $$
Note that $\lim_{\gamma\to\infty}c^-_{\gamma,d}=1$.
From \eqref{finallower} and the above remarks, we have now shown that
\begin{equation}\label{largegamma}
G^D_{\text{Neu}}(x,y)\ge\frac1{(d-2)\omega_d}c^-_{\gamma,d}|y|^{2-d}, \ \text{for}\ x\in\partial D, |y|\ge\gamma R, \gamma\ge \gamma_0.
\end{equation}
This is the lower bound in \eqref{key} except that one has here $\gamma\ge\gamma_0$ instead of $\gamma>1$.
To extend \eqref{largegamma} to $\gamma>1$, one notes
that
since $G^R_{\text{Dir}}(z,y)$ is positive for $z,y\in R^d-\bar B_R$,
 one has trivially  $\inf_{\gamma R\le |y|\le \gamma_0 R, |z|=\rho R} G^R_{\text{Dir}}(z,y)>0$, for
any choice of $\rho, \gamma$ satisfying
 $1<\rho<\gamma<\gamma_0$,
 and in fact by scaling, the left hand side is independent of $R$.
Using this along with the fact that
$G^D_{\text{Neu}}(x,A)\ge\left(\inf_{z\in \partial B_{\rho R}}\int_AG^R_{\text{Dir}}(z,y)dy\right)\frac1{1-\rho^{2-d}}$, which follows from
\eqref{finallower}, we can define $c^-_{\gamma,d}>0$ for all $\gamma>1$ so that \eqref{largegamma} holds.

\end{proof}

\section{Proof of Theorem \ref{LHE-gen}}\label{sec:gen}
In the present case,  \eqref{formula}  still holds, where  $G^D_{\text{Neu}}$
is now the Green's function
for the  operator $L$ on $R^d-\bar D$ with
 the co-normal  boundary condition at $\partial D$.
 Since the operator $L$ is symmetric with respect to the weight $e^Q$, similar to \eqref{sym}
we have
\begin{equation}\label{Qsym}
 e^{Q(x)}G^D_{\text{Neu}}(x,y)=e^{Q(y)}G^D_{\text{Neu}}(y,x).
\end{equation}
To see this, note that since Green's function $G^D_{\text{Neu}}(x,y)$ satisfies $LG^D_{\text{Neu}}(\cdot, y)=-\delta_y$  with the homogeneous co-normal boundary condition
at $\partial D$,
it follows that  for any compactly supported sufficiently smooth function $f$ defined
on $R^d-D$,
the function $v(x)=\int_{R^d-D}G^D_{\text{Neu}}(x,y)f(y)dy$ is the minimal positive solution of
$Lv=-f$ in $R^d-D$,  with the homogeneous co-normal boundary condition at $\partial D$. Now
on the one hand, since $L$ is in the form \eqref{Lsym}, it follows that $v$ solves
$\nabla\cdot e^Qa\nabla v=-e^Qf$ with the homogeneous co-normal boundary condition,
but on the other hand, by the same reasoning as above,
we have $v(x)=\int_{R^d-D} G_{\text{sym}}(x,y)e^Q(y)f(y)dy$,
where $G_{\text{sym}}$ is the Green's function for the operator
$\nabla\cdot e^Qa\nabla$ on $R^d-D$  with the
co-normal boundary condition at $\partial D$.
Since this holds for all nice $f$, we conclude
that $G^D_{\text{Neu}}(x,y)=G_{\text{sym}}(x,y)e^{Q(y)}$.
Since $G_{\text{sym}}$ is the Green's function of an operator that is symmetric with respect to Lebesgue measure, it follows
as in \eqref{sym} that $G_{\text{sym}}(x,y)=G_{\text{sym}}(y,x)$.
Now \eqref{Qsym} follows from this.

In light of \eqref{formula} and \eqref{Qsym},
 to prove the theorem it suffices to show that for $R'>R$, one has
\begin{equation}\label{need}
\begin{aligned}
\frac{\min_{|z|=R'}G^R_{\text{Dir}}(z,x)}{1-\min_{|z|=R'}V(z)}
\le G^D_{\text{Neu}}(y,x)\le
\frac{\max_{|z|=R'}G^R_{\text{Dir}}(z,x)}{1-\max_{|z|=R'}V(z)},\   \text{for}\ |x|>R' \ \text{and}\ y\in\partial D.
\end{aligned}
\end{equation}

As in the previous section, we consider  the probabilistic representation of the Green's function.
 Let $X(t)$ be the diffusion process 
in $R^d-D$, co-normally reflected at $\partial D$, and corresponding to the operator $L$ 
\cite{SV}.
Let $P_y$ denote probabilities and let $E_y$ denote the corresponding  expectations for the diffusion starting from $y\in R^d-D$.
 For $y\in R^d-D$, define the expected occupation measure
 by $\mu_y^D(A)=E_y\int_0^\infty1_A(X(t))dt$, for Borel sets $A\subset R^d-D$.
 By assumption, the process $X(t)$ is transient,
 (that is, $P_y(\lim_{t\to\infty}|X(t)|=\infty)=1$), and from this one can  show that  $\mu_y^D(A)<\infty$ for all bounded $A$.
 The measure $\mu_y^D(dx)$ possesses  a density
and the density is given by $G^D_{\text{Neu}}(y,x)$.
From now on  we will write $G^D_{\text{Neu}}(y,A)\equiv \mu_y^D(A)$.

Define
  $\tau_1=\inf\{t\ge0:|X(t)|=R_1\}$,  and then by induction define
$\tau_{2n}=\inf\{t>\tau_{2n-1}:|X(t)|=R\}$ and $\tau_{2n+1}=\inf\{t>\tau_{2n}: |X(t)|= R_1\}$.
Similar to  the proof of Theorem \ref{LHE-Lap} (see \eqref{circuits}, \eqref{avg} and the line between \eqref{free} and \eqref{final}),
we have for $y\in\partial D$ and $A\subset R^d-\bar B_{R_1}$,
\begin{equation}\label{long}
\begin{aligned}
&G^D_{\text{Neu}}(y,A)=E_y\int_0^\infty1_A(X(t))dt=\sum_{n=1}^\infty E_y\int_{\tau_{2n-1}}^{\tau_{2n}}1_A(X(t))dt=\\
&\sum_{n=1}^\infty E_y(\int_{\tau_{2n-1}}^{\tau_{2n}}1_A(X(t))dt|\tau_{2n-1}<\infty)P_y(\tau_{2n-1}<\infty)=\\
&\sum_{n=1}^\infty\int_{\partial B_{R_1}}\left(E_z\int_0^{\tau_2}1_A(X(t))dt\right)\nu_{n,y}(dz)P_y(\tau_{2n-1}<\infty),
\end{aligned}
\end{equation}
where $\nu_{n,y}$ is some  probability measure on $\partial B_{R_1}$.
By the strong Markov property,
\begin{equation}\label{twosidedhitting}
(\min_{|z|=R'}V(z))^{n-1}\le P_y(\tau_{2n-1}<\infty)\le(\max_{|z|=R'}V(z))^{n-1},
\end{equation}
for $y\in\partial D$, where $V$ is as in \eqref{Vprob}.
Now \eqref{need} follows from \eqref{long}, \eqref{twosidedhitting} and the fact that
$E_z\int_0^{\tau_2}1_A(X(t))dt=G_{\text{Dir}}^R(z,A)$.

\noindent \bf Acknowledgment.\rm\ The author thanks Gilbert Weinstein for bringing some of the aspects of  this problem
 to his attention.


\begin{thebibliography}{99}
\bibitem{GT} Gilbarg, D. and  Trudinger, N.,
 Elliptic Partial Differential Equations of Second Order, second edition, Springer-Verlag, 1983

\bibitem{J}
John, F.,  Partial Differential Equations,  Third Edition,  Applied Mathematical Sciences, 1,  Springer-Verlag, 1978.

\bibitem{P} Pinsky R., Positive Harmonic Functions and Diffusion, Cambridge University Press, 1995.

\bibitem{SV} Stroock, D. and Varadhan, S.R.S., \emph{Diffusion processes with boundary conditions}, Comm. Pure Appl. Math. 24 1971 147-225. 



\end{thebibliography}
\end{document}